\documentclass[12pt,a4paper, final]{article}

\usepackage[latin1]{inputenc}

\usepackage[T1]{fontenc}
\usepackage[english]{babel}
\usepackage{amsmath,amssymb,amscd,latexsym}

\input{diagrams}

\usepackage{theorem}
\theorembodyfont{\slshape}
\newtheorem{Thm}{Theorem}

\theoremstyle{plain}

\newtheorem{Prop}{Proposition}
\newtheorem{Lem}{Lemma}

\theorembodyfont{\rmfamily}

\newcommand{\Z}{\mathbb{Z}}
\newcommand{\R}{\mathbb{R}}
\newcommand{\F}{\mathbb{F}}

\newcommand{\qed}{\hfill $\square$}
\newcommand{\Pres}[2]{\left\langle{#1}\ \big\vert\ {#2}\right\rangle}

\renewcommand{\k}{\mathfrak k}

\newcommand{\m}{\mathfrak m}
\newcommand{\n}{\mathfrak n}

\title{Convergence of Baumslag-Solitar groups}
\author{Yves STALDER\footnote{Supported by the Swiss national Science Foundation, grant 20-101469}}
\date{January 17, 2005}

\begin{document}

\selectlanguage{english}

\maketitle

\begin{abstract}
We study convergent sequences of Baumslag-Solitar groups in the space of marked groups. We prove that
$BS(\m,\n) \to \F_2$ for $|\m|,|\n| \to \infty$ and $BS(1,\n) \to \Z\wr\Z$ for $|\n| \to \infty$. For $\m$ fixed,
$|\m| \geqslant 2$, we show that the sequence $(BS(\m,\n))_{\n}$ is not convergent and characterize many
convergent subsequences. Moreover if $X_\m$ is the set of $BS(\m,\n)$'s for $\n$ relatively prime to $\m$ and $|\n|
\geqslant 2$, then the map $BS(\m,\n) \mapsto \n$ extends continuously on $\overline{X_\m}$ to a surjection onto
invertible $\m$-adic integers.
\end{abstract}

\section{Introduction}\label{Intro}

Let ${\cal G}_2$ be the space of finitely generated marked groups on two generators (see Section \ref{Def} for
definition) and let $\F_2 = \Pres{a,b}{\varnothing}$ be the free group on two generators. Baumslag-Solitar groups
are defined by the presentations
$$
BS(\m,\n)= \Pres{a,b}{ab^{\m}a^{-1} = b^\n}
$$
for $\m,\n \in \Z^* = \Z \setminus \{0\}$. The purpose of the present paper is to understand how Baumslag-Solitar
groups are distributed in ${\cal G}_2$. More precisely, we determine convergent sequences and in some cases we
are able to give the limit group. In the following results, we mark $\F_2$ and $BS(\m,\n)$ by $\{a,b\}$.

\begin{Thm}\label{ToFree}
For $|\m|,|\n| \to \infty$, we have $BS(\m,\n) \to \F_2$.
\end{Thm}

In particular, the property of being Hopfian is not open in ${\cal G}_2$ since $BS(2^k,3^k)$ is known to be
non Hopfian for all $k\geqslant 1$ (see \cite{LynSch}, Chapter IV, Theorem 4.9.) while $\F_2$ is Hopfian.
Theorem \ref{ToFree} is not so surprising because the length of the relator appearing in the presentation of
$BS(\m,\n)$ tends to $\infty$  when $|\m|,|\n| \to \infty$. However, this relator is not the shortest relation in
the group for many values of $(\m,\n)$. To prove Theorem \ref{ToFree}, we give a lower bound for the length of
shortest relations in the more general setting of HNN-extensions (see Section \ref{LowBd}).

We now fix the parameter $\m$. In the case $\m = \pm 1$, we show that the sequence $(BS(\pm 1,\n))_{\n\in\Z}$
is convergent and we can identify the limit.

\begin{Thm}\label{ToWreath}
Let the wreath product $\Z \wr \Z = \Z \ltimes \oplus_{i\in\Z} \Z$ be marked by the elements $(1,0)$ and
$(0,e_0)$ where $e_0 \in \oplus_{i\in\Z} \Z$ is the Dirac mass at $0$. Then $BS(\pm 1,\n) \to \Z \wr \Z$ when
$|\n| \to \infty$.
\end{Thm}

This illustrates the fact that a limit of metabelian groups is metabelian, see \cite{ChaGui}, Section 2.6.

In the case $|\m| \geqslant 2$, we show that the sequence $(BS(\m,\n))_{\n\in\Z}$ is not convergent in ${\cal
G}_2$. As ${\cal G}_2$ is compact, it has convergent subsequences. The next result we state in this
introduction, among subsequences, characterizes many convergent ones. However we don't actually know what the
limits are. Notice that the result also holds for $\m= \pm 1$, even if it is in this case weaker than Theorem
\ref{ToWreath}.

\begin{Thm}\label{ConvSub}
Let $\m\in\Z^*$ and let $(\k_n)_n$ be a sequence of integers relatively prime to $\m$. The sequence $(BS(\m,\k_n))_n$ is
convergent in ${\cal G}_2$ if and only if one (and only one) of the following assertions holds:
\begin{enumerate}
  \item[(a)] $(\k_n)_n$ is eventually constant;
  \item[(b)] $|\k_n| \to \infty$ and $(\k_n)_n$ is eventually constant
  modulo $\m^h$ for all $h\geqslant 1$.
\end{enumerate}
\end{Thm}

Note that for $|\m| \geqslant 2$, condition (b) precisely means that $|\k_n| \to \infty$ and $(\k_n)_n$ is
convergent in $\Z_\m$, the ring of $\m$-adic integers. The link between Baumslag-Solitar groups and $\m$-adic
integers can be made more precise. We define $X_\m$ to be the set of $BS(\m,\n)$'s, for $\n$ relatively prime to
$\m$ and $|\n| \geqslant 2$ and we denote by $\Z_\m^{\times}$ the set of invertible elements of $\Z_\m$.

\begin{Thm}\label{UnifCont}
For $|\m| \geqslant 2$, the map $\Psi : X_\m \to \Z_\m^\times \ ; \ BS(\m,\n) \mapsto \n$ extends
continuously to $\overline{X_\m}$. The extension is surjective, but not injective.
\end{Thm}

An immediate corollary of Theorem \ref{ConvSub} or Theorem \ref{UnifCont} is that (for $|\m|\geqslant 2$) the
sequence $(BS(\m,\n))_\n$ admits uncountably many accumulation points, namely at least one for each invertible
$\m$-adic integer.

We end this introduction by a remark on markings of Baumslag-Solitar groups. In this paper we always mark the
group $BS(\m,\n)$ by the generators coming from its canonical presentation given above. Nevertheless, it is
also an interesting approach to consider different markings on $BS(\m,\n)$. For instance, take $\m$ and $\n$
greater than $2$ and relatively prime, so that $\Gamma = BS(\m,\n)$ is non-Hopfian, the epimorphism
$\phi:\Gamma \to \Gamma$ given by $a \mapsto a$ and $b \mapsto b^\m$ being non-injective (see again
\cite{LynSch}, Chapter IV, Theorem 4.9.). In \cite{ABLRSV}, the authors consider the sequence of groups
$\Gamma_n = \Gamma/\ker(\phi^n)$ (marked by $a$ and $b$). They show that the sequence $(\Gamma_n)_n$
converges to an amenable group while, being all isomorphic to $\Gamma$ as groups, the $\Gamma_n$'s are not
amenable. This allow them to prove that $\Gamma$ is non-amenable, but not uniformly (cf. Proposition 13.3)
and also shows that the property of being amenable is not open in ${\cal G}_2$.

\paragraph{Structure of the article.} In Section \ref{Def}, we give the necessary preliminaries. In Section
\ref{LowBd}, we estimate the length of the shortest relations in a HNN-extension and prove Theorem
\ref{ToFree}. Section \ref{m1} is devoted to the case $\m =1$, namely the proof of Theorem \ref{ToWreath},
and Section \ref{1pBS} to other one-parameter families of Baumslag-Solitar groups, i.e. the proof of Theorem
\ref{ConvSub}. Finally, we prove Theorem \ref{UnifCont} in Section \ref{Topol}, linking Baumslag-Solitar
groups and $\m$-adic integers.

\paragraph{Acknowledgements.} I would like to thank Luc Guyot and Alain Valette for their useful comments and hints,
and Indira Chatterji, Pierre de la Harpe, Nicolas Monod and Alain Robert for comments on previous versions.

\section{Preliminaries}\label{Def}

In this Section, we collect some definitions and material which are needed in the rest of the paper. The
reader who is familiar with the notions of $\m$-adic integers, HNN-extensions and
topology on the space of marked groups can skip this Section.

\paragraph{The ring of $\m$-adic integers.}
For $\m \in \Z$, $|\m| \geqslant 2$ we define $\Z_\m$ to be the completion of $\Z$
with respect to the ultrametric distance given by the following "absolute value":
$$
|a \cdot \m^{k}|_\m := \left( \frac{1}{|\m|} \right)^{k} \ \text{ for } a \text{ not a multiple of } \m \text
{ and } k\geqslant 0 \ .
$$
(Note that it is not multiplicative in general, but one has $|-x|_\m = |x|_\m$ which is sufficient to induce a distance.)
The space $\Z_\m$ has a ring structure obtained by continuous extensions of the ring laws for $\Z$ and we call it
the \emph{ring of $\m$-adic integers}. The symbol $\Z_\m^\times$ denotes the invertible elements of $\Z_\m$. The
distance induced by $|.|_\m$ is called the \emph{$\m$-adic ultrametric distance}, since it satisfies the ultrametric inequality
$$
|x-z|_\m \leqslant \text{max} \big( |x-y|_\m,|y-z|_\m \big) \ .
$$
As a topological ring, $\Z_\m$ is the projective limit of the system
$$
\ldots \to \Z/\m^h\Z \to \Z/\m^{h-1}\Z \to \ldots \to \Z/\m^2\Z \to \Z/\m\Z
$$
where the arrows are the canonical (surjective) homomorphisms. This shows that $\Z_\m$ is compact and it
is coherent with this characterization to set $\Z_{\m} = \{0\}$ for $\m = \pm 1$. It becomes also
nearly obvious that the group of invertible elements of $\Z_\m$ is given by
$$
\Z_\m^\times = \Z_\m \setminus (p_1 \Z_\m \cup \ldots \cup p_k \Z_\m)
$$
where $p_1, \ldots, p_k$ are the prime factors of $\m$.

To conclude this short summary about $\m$-adic integers, let us notice that, for $|\m| >1$
and $\m$ not prime, the ring $\Z_\m$ has zero divisors.

\paragraph{Marked groups and their topology.} Introductory expositions of these topics can be found in
\cite{Cha} or \cite{ChaGui}. We only recall some basics and what we need in following sections.

The free group on $k$ generators will be denoted $\F_k$, or $F_S$ (with $S = (s_1, \ldots, s_k)$)
if we want to precise the names of (canonical) generating elements.
A \emph{marked group on $k$ generators} is a pair $(\Gamma,S)$ where $\Gamma$ is a group and
$S = (s_1, \ldots, s_k)$ is a family which generates $\Gamma$. A marked group $(\Gamma,S)$ comes always
with a canonical epimorphism $\phi: \F_S \to \Gamma$, giving an isomorphism of marked groups between
$\F_S/ \ker \phi$ and $\Gamma$. Hence a class of marked groups can always be represented by a quotient of
$\F_S$. In particular if a group is given by a presentation, this defines a marking on it. The nontrivial elements of
${\cal R} := \ker \phi$ are called \emph{relations} of $(\Gamma,S)$.

Let $w = x_1^{\varepsilon_1} \cdots x_n^{\varepsilon_n}$ be a reduced word in $\F_S$ (with $x_i \in S$ and
$\varepsilon_i \in \{ \pm 1 \}$). The integer $n$ is called the \emph{length} of $w$ and denoted $\ell(w)$.
The length of the shortest relation(s) of $\Gamma$ will be denoted $g_\Gamma$, for we observe it is the girth
of the Cayley graph of $\Gamma$ (with respect to the generating set $S$). In case ${\cal R} = \varnothing$,
we set $g_\Gamma = +\infty$.

If $(\Gamma,S)$ is a marked group on $k$ generators, and $\gamma \in \Gamma$ the \emph{length} of $\gamma$ is
\begin{eqnarray*}
\ell_\Gamma(\gamma) & := & \min\{ n: \gamma = s_1 \cdots s_n \text{ with } s_i \in S \sqcup S^{-1} \} \\
                    &  = & \min\{ \ell(w) : w \in \F_S, \ \phi(w) = \gamma  \} \ .
\end{eqnarray*}

Let ${\cal G}_k$ be the set of marked groups on $k$ generators (up to marked isomorphism). Let us recall that
the topology on ${\cal G}_k$ comes from the following ultrametric: for $(\Gamma_1, S_1) \neq (\Gamma_2, S_2)
\in {\cal G}_k$ we set $d \big( (\Gamma_1, S_1), (\Gamma_2, S_2) \big) := e^{-\lambda}$ where $\lambda$ is
the length of a shortest element of $\F_k$ which vanishes in one group and not in the other one. But what the
reader has to keep in mind is the following characterization of convergent sequences.

\begin{Prop}
Let $(\Gamma_n, S_n)$ be a sequence of marked groups (on $k$ generators). The following are equivalent:
\begin{enumerate}
  \item[(i)] $(\Gamma_n, S_n)$ is convergent in ${\cal G}_k$;
  \item[(ii)] for all $w \in \F_k$ we have either $w = 1$ in $\Gamma_n$ for $n$ large enough,
  or $w \neq 1$ in $\Gamma_n$ for $n$ large enough.
\end{enumerate}
\end{Prop}

\textbf{Proof.} (i)$\Rightarrow$(ii): Set $(\Gamma,S) = \underset{n\to \infty}{\lim} (\Gamma_n,S_n)$ and take
$w \in \F_k$. For $n$ sufficiently large we have $d\big( (\Gamma,S),(\Gamma_n,S_n) \big) < e^{- \ell(w)}$,
which implies that we have $w=1$ in $\Gamma_n$ if and only if $w=1$ in $\Gamma$.

(ii)$\Rightarrow$(i): Set $N = \{ w \in \F_k : w=1 \text{ in } \Gamma_n \text{ for } n \text{ large enough }
\}$, $\Gamma = \F_k/N$, and fix $r \geqslant 1$. For $n$ large enough, $\Gamma_n$ and $\Gamma$ have the same
relations up to length $r$ (for the balls in $\F_k$ are finite) and hence $d(\Gamma,\Gamma_n) < e^{- r}$ (we
drop the markings since they are obvious). This implies $\Gamma_n \underset{n \to \infty}{\to} \Gamma$. \qed

\paragraph{HNN-extensions and Baumslag-Solitar groups.} Suppose now that $(H,S)$ is a marked group on $k$
generators, and that $\phi: A \to B$ is an isomorphism between subgroups of $H$. The \emph{HNN extension}
of $H$ with respect to $A$, $B$ and $\phi$ is given by
$$
HNN(H,A,B,\phi) := \frac{H * \langle t \rangle}{\cal N} \ .
$$
where ${\cal N}$ is the normal subgroup generated by the $t^{-1}at \phi(a)^{-1}$ for $ a \in A$.
Unless specified otherwise, we always mark a HNN-extension by $S$ and $t$. An element $\gamma \in HNN(H,A,B,\phi)$
can always be written
\begin{equation}\label{decHNN}
  \gamma = h_0 t^{\varepsilon_1} h_1 \cdots t^{\varepsilon_n} h_n \text{ with } n \geqslant 0, \
  \varepsilon_i \in \{\pm 1\}, \ h_i \in H\ .
\end{equation}
The decomposition of $\gamma$ in (\ref{decHNN}) is called \emph{reduced} if no subword of type
$t^{-1}at$ (with $a \in A$) or $tbt^{-1}$ (with $b \in B$) appears. We recall the following result, which
is called \emph{Britton's Lemma}.
\begin{Lem}\label{Britton}
(\cite{LynSch}, Chapter IV.2.) Let $\gamma \in HNN(H,A,B,\phi)$ and write as in (\ref{decHNN}) $\gamma =
h_0 t^{\varepsilon_1} h_1 \cdots t^{\varepsilon_n} h_n$. If $n \geqslant 1$ and if the decomposition is
reduced, then $\gamma \neq 1$ in $HNN(H,A,B,\phi)$.
\end{Lem}

This shows in particular that the integer $n$ appearing in a reduced decomposition is uniquely determined by $\gamma$.

Let us finally recall that \emph{Baumslag-Solitar groups} are defined by the presentations $BS(\m,\n)= \Pres{a,b}{ab^{\m}a^{-1} =
b^\n}$ for $\m,\n \in \Z^*$. Setting $\phi(\n k) = \m k$, we have $BS(\m,\n) = HNN(\Z, \n\Z, \m\Z, \phi)$.

\section{Shortest relations in a HNN-extension and convergence of Baumslag-Solitar groups}\label{LowBd}

Let $(H,S)$ be a marked group and $\Gamma = HNN(H,A,B,\phi)$. In this Section we give a lower estimate for
$g_\Gamma$. As a higher estimate, we obviously get $g_\Gamma \leqslant g_H$, because a shortest relation
in $H$ is also a relation in $\Gamma$. Let us define:
\begin{eqnarray*}
    \alpha & := & \min\big\{ \ell_H(a) \ : \ a \in A \setminus \{1\} \big\} \ ;   \\
    \beta  & := & \min\big\{ \ell_H(b) \ : \ b \in B \setminus \{1\} \big\} \ .
\end{eqnarray*}
\begin{Thm}\label{Low}
Let $(H,S)$, $\Gamma$, $\alpha$ and $\beta$ be defined as above. Then
we have
$$
\min\{ g_H, \ \alpha+\beta+2, \ 2\alpha+6, \ 2\beta+6 \} \leqslant g_\Gamma \leqslant g_H \ .
$$
\end{Thm}

As the case of Baumslag-Solitar groups (treated below) will show, the lower bound given in
Theorem \ref{Low} is in fact sharp. This sharpness is the principal interest of this Theorem, because the
estimate $\min\{ g_H, \ \alpha, \ \beta \} \leqslant g_\Gamma$, which follows easily from Lemma \ref{rel}
and Britton's Lemma, suffices to prove Theorem \ref{ToFree} (replace Proposition \ref{BSgirth} by the estimate
$g_{BS(\m,\n)} \geqslant \min(\m,\n)$). I would like to thank the referee for having pointed
this fact to me.

Before proving Theorem \ref{Low}, let us begin with a simple observation.
\begin{Lem}\label{rel}
Let $(H,S)$ and $\Gamma$ be defined as above and let $r$ be a relation of $\Gamma$ contained in $\F_S$.
Then $r$ is a relation of $H$. In particular, $\ell(r) \geqslant g_H$.
\end{Lem}

\textbf{Proof.} Since $r=1$ in $\Gamma$ and since the canonical map $H \to \Gamma$ is injective, we get $r=1$ in $H$.
Hence the first assertion. The second one follows by definition of $g_H$. \qed

\smallskip

\textbf{Proof of Theorem \ref{Low}.} The second inequality has already been discussed. To establish the first one, let us take a
relation $r$ of $\Gamma$ and show that $\ell(r) \geqslant m$, where we set $m := \min\{ g_H, \ \alpha+\beta+2,
\ 2\alpha+6, \ 2\beta+6 \}$.

Write $r = h_0 t^{\varepsilon_1} h_1 \cdots t^{\varepsilon_n} h_n$ with $\varepsilon_i \in \{\pm 1\}$, $h_i
\in \F_S$ and $h_i \neq 1$ if $\varepsilon_i = - \varepsilon_{i+1}$. Up to replacement by a (shorter)
conjugate, we may assume that $r$ is cyclically reduced. If $n \neq 0$, we may also assume that $h_0 = 1$.
Since $r=1$ in $\Gamma$, one clearly has $\sum_{i=1}^{n} \varepsilon_i = 0$. In particular, $n$ is even. Let
us distinguish several cases and show $\ell(r) \geqslant m$ in each one:

\textbf{Case $n=0$:} We get $r = h_0 \in \F_S$. Thus $\ell(r) \geqslant g_H \geqslant m$ by Lemma \ref{rel}.

\textbf{Case $n=2$:} One gets $r = t^{\varepsilon} h_1 t^{-\varepsilon} h_2$. If we look at $r$ in $\Gamma$, we have
$r=1$ and thus $h_1 \in A$ (if $\varepsilon = -1$) or $h_1 \in B$ (if $\varepsilon =1$) by Britton's Lemma. Suppose
$\varepsilon = -1$ (in case $\varepsilon = 1$ the proof is similar and left to the reader). Looking at $h_1$ in $\F_S$,
there are two possibilities (remember that we assumed $h_1 \neq 1$).
\begin{itemize}
  \item If $h_1 = 1$ in $\Gamma$, Lemma \ref{rel} implies $\ell(r) \geqslant \ell(h_1) \geqslant g_H \geqslant m$.
  \item If $h_1 \neq 1$ in $\Gamma$, then $\ell(h_1) \geqslant \alpha$. On the other hand $h_2^{-1} = t^{-1} h_1 t \in
  B$; thus $\ell(h_2) \geqslant \beta$ and $\ell(r) \geqslant \alpha + \beta + 2 \geqslant m$.
\end{itemize}

\textbf{Case $n \geqslant 4$:} We have $r = 1$ in $\Gamma$. By Britton's Lemma, there is an index $i$ such
that either $\varepsilon_i = -1$, $\varepsilon_{i+1} = 1$ and $h_i \in A$, or $\varepsilon_i = 1$,
$\varepsilon_{i+1} = -1$ and $h_i \in B$. Since cyclic conjugations preserve length, we may assume $i=1$, so that
$r = t^{\varepsilon_1} h_1 t^{-\varepsilon_1} h_2 t^{\varepsilon_3} h_3 \cdots t^{\varepsilon_n} h_n$. Let us
moreover assume $\varepsilon_1 = -1$ (again the case $\varepsilon_1 = 1$ is similar and left to the reader).
Set $r' = w h_2 t^{\varepsilon_3} h_3 \cdots t^{\varepsilon_n} h_n$ where $w \in F_X$ is such that $w =
t^{-1}h_1t$ in $\Gamma$ (in fact this element is in $B$). Applying Britton's Lemma to $r'$, one sees there
exists an index $j \geqslant 3$ such that either $\varepsilon_{j} = -1 = -\varepsilon_{j+1}$ and $h_j \in A$,
or $\varepsilon_j = 1 = -\varepsilon_{j+1}$ and $h_j \in B$.  There are three possibilities:
\begin{itemize}
  \item If $h_1 = 1$ or $h_j = 1$ in $\Gamma$, Lemma \ref{rel} implies $\ell(r) \geqslant g_H \geqslant m$ as
  above.
  \item If $h_1 \neq 1$ in $\Gamma$, $h_j \neq 1$ in $\Gamma$ and $\varepsilon_j = 1$, then $\ell(h_1) \geqslant \alpha$
  and $\ell(h_j) \geqslant \beta$. Thus $\ell(r) \geqslant \alpha + \beta + 4 > m$.
  \item If $h_1 \neq 1$ in $\Gamma$, $h_j \neq 1$ in $\Gamma$ and $\varepsilon_j = -1$, then we write $r =
  t^{-1}h_1t w_1 t^{-1}h_jt w_2$ with $\ell(h_1) \geqslant \alpha$ and $\ell(h_j) \geqslant \alpha$. The
  subwords $w_1, w_2$ are not empty because $r$ is cyclically reduced. Thus $\ell(r) \geqslant 2 \alpha + 6
  \geqslant m$ (Remark that $2\alpha + 6$ would be replaced by $2 \beta +6$ in the case $\varepsilon_1 = 1,
  \varepsilon_j = 1$).
\end{itemize}
The proof is complete. \qed

\smallskip

We now turn to prove that for Baumslag-Solitar groups, the lower bound coming
from Theorem \ref{Low} is in fact the length of the shortest relation. More precisely we have the following
statement:
\begin{Prop}\label{BSgirth}
Let $\m,\n \in \Z^*$. We have
$$
g_{BS(\m,\n)} = \min\big\{ |\m|+|\n|+2, \ 2|\m|+6, \ 2|\n|+6 \big\} \ .
$$
\end{Prop}

\textbf{Proof.} Set $m := \min\{ |\m|+|\n|+2, \ 2|\m|+6, \ 2|\n|+6 \}$ and $\Gamma = BS(\m,\n)$. We have $g_\Z =
+\infty$, $\alpha = |\n|$ and $\beta = |\m|$. Thus, Theorem \ref{Low} implies $g_\Gamma \geqslant m$. To prove
that $g_\Gamma \leqslant m$, we produce relations of length $|\m|+|\n|+2$, $2|\m|+6$, and $2|\n|+6$. Namely:
\begin{itemize}
    \item $ab^\m a^{-1}b^{-\n}$ has length $|\m|+|\n|+2$;
    \item $ab^\m a^{-1}b ab^{-\m}a^{-1}b^{-1}$ has length $2|\m|+6$;
    \item $a^{-1}b^\n ab a^{-1}b^{-\n}ab^{-1}$ has length $2|\n|+6$. \qed
\end{itemize}

Theorem \ref{ToFree} of introduction is now a consequence of Proposition \ref{BSgirth}, since
a sequence of groups $\Gamma_n = \Pres{a,b}{{\cal R}_n}$ converges to the free group on two
generators (marked by its canonical basis) if and only if $g_{\Gamma_n}$ tends to $\infty$.
\section{Limit of Solvable Baumslag-Solitar groups}\label{m1}

This section is entirely devoted to the proof of Theorem \ref{ToWreath}. We notice first that
$BS(1,\n) = BS(-1,-\n)$ as marked groups. Thus we
may assume $\m = 1$. Hence, we let $\Gamma_\n = BS(1,\n)$ and $\Gamma = \Z \wr \Z$. In $\Z \wr \Z$, let us set
$a= (1,0)$ and $b = (0,e_0)$. We have to show that for all $w \in \F_2$:
\begin{enumerate}
  \item[(1)] if $w=1$ in $\Z \wr \Z$, then $w=1$ in $BS(1,\n)$ for $|\n|$ large enough;
  \item[(2)] if $w\neq 1$ in $\Z \wr \Z$, then $w\neq 1$ in $BS(1,\n)$ for $|\n|$ large enough.
\end{enumerate}

Let $w \in \F_2$. One can write $w = a^\alpha a^{\alpha_1} b^{\beta_1} a^{-\alpha_1} \cdots
a^{\alpha_k} b^{\beta_k} a^{-\alpha_k}$. The image of $w$ in $\Gamma$ is $(\alpha, \sum_{i=1}^k \beta_i
e_{\alpha_i})$, where $e_j \in \oplus_{h \in \Z} \Z$ is the Dirac mass at $j$.
Let $m = \min_{1 \leqslant i \leqslant k} \alpha_i$. In $\Gamma_\n = BS(1,\n)$, we have
\begin{eqnarray*}
    w & = & a^\alpha a^m a^{\alpha_1-m} b^{\beta_1} a^{m-\alpha_1} \cdots a^{\alpha_k-m} b^{\beta_k} a^{m-\alpha_k} a^{-m}  \\
      & = & a^\alpha a^m b^{\beta_1 \n^{\alpha_1-m}} \cdots b^{\beta_k \n^{\alpha_k-m}} a^{-m}  \\
      & = & a^\alpha a^m b^{\sum_{h \in \Z} \left( \sum_{\alpha_i = h} \beta_i \right) \n^{h-m}} a^{-m}
\end{eqnarray*}

(1) As $w \underset{\Gamma}{=} 1$, we have $\alpha = 0$ and $\forall h \in \Z$, $\sum_{\alpha_i = h}
\beta_i = 0$. Hence
$$
w \underset{\Gamma_\n}{=} a^0 a^m b^{\sum_{h \in \Z} 0 \cdot \n^{h-m}} a^{-m} = 1 \ \forall \n\in \Z^* \ .
$$

(2) As $w \underset{\Gamma}{\neq} 1$, either $\alpha \neq 0$ or $\exists h \in \Z$ such that $\sum_{\alpha_i = h}
\beta_i \neq 0$. The image of $w$ by the morphism $\Gamma_\n \to \Z$ given by $a \mapsto
1, b \mapsto 0$ is $\alpha$. Hence, if $\alpha \neq 0$, then $w \underset{\Gamma_\n}{\neq} 1 \ \forall \n$. If
$\alpha = 0$, we set $h_0$ to be the maximal value of $h$ such that $\sum_{\alpha_i = h} \beta_i \neq 0$. For
$|\n|$ large enough, we have
$$
\left| \sum\limits_{\alpha_i = h_0} \beta_i \n^{h_0 - m} \right| > \left| \sum\limits_{h<h_0} \sum\limits_{\alpha_i = h}
\beta_i \n^{h-m} \right| \ .
$$
For those values of $\n$, we get
$$
w \underset{\Gamma_\n}{=} a^m b^{\left( \sum_{\alpha_i = h_0} \beta_i \right) \n^{h_0 - m} +
\sum_{h < h_0} \left( \sum_{\alpha_i = h} \beta_i \right) \n^{h-m}} a^{-m} \underset{\Gamma_\n}{\neq} 1 \ .
$$
The proof is complete. \qed

\section{General one-parameter families of Baumslag-Solitar groups}\label{1pBS}

We now treat the case $|\m| \geqslant 2$. More precisely, we begin the proof of Theorem \ref{ConvSub}. We also
have $BS(\m,\n) = BS(-\m,-\n)$ as marked groups. This equality will allow us to assume $\m > 0$ in following
proofs. We begin with a Lemma which already shows that the sequence
$(BS(\m,\n))_\n$ is not itself convergent.

\begin{Lem}\label{Congruence}
    Let $\m, \n \in \Z^*$, $d= \text{gcd}(\m, \n)$. We write $\m = d \m_1$, $\n = d \n_1$. Let $k\in \Z$, $h\geqslant 1$ and
    $$
        w = a^{h+1} b^{\m} a^{-1} b^{-k} a^{-h} b a^{h+1} b^{-\m} a^{-1} b^k a^{-h} b^{-1} \ .
    $$
    If $|\n| \geqslant 2$, we have $w = 1$ in $BS(\m, \n)$ if and only if $\n \equiv k \
    (\text{mod } \m_1^h d)$.
\end{Lem}

The congruence modulo $\m_1^h d$ (instead of $\m^h$) is the reason for the hypothesis "$\n$ relatively prime
to $\m$" appearing in Theorem \ref{ConvSub}.

\smallskip

\textbf{Proof.} Let $\Gamma_{\n} = BS(\m, \n)$. We have
$$
    w \underset{\Gamma_{\n}}{=} a^h b^{\n - k} a^{-h} b a^{h} b^{k-\n} a^{-h} b^{-1} \ .
$$
Let us now distinguish three cases:

\textbf{Case 1:} $\n \not\equiv k \ (\text{mod } \m)$. \\
We have $w \underset{\Gamma_{\n}}{\neq} 1$ by Britton's Lemma, since $|\n| \geqslant 2$.

\textbf{Case 2:} $\n \not\equiv k \ (\text{mod } \m_1^h d)$, but $\n \equiv k \ (\text{mod } \m)$. \\
We write $\n - k = \ell \m_1^g d$ with $g<h$ and $\ell$ not a multiple of $\m_1$. Hence $\ell \n_1^g d$ is not a multiple
of $\m = \m_1 d$, for $\m_1$ is relatively prime to $\n_1$. We have
$$
w \underset{\Gamma_{\n}}{=} a^{h-g} b^{\ell \n_1^g d} a^{g-h} b a^{h-g} b^{-\ell \n_1^g d} a^{g-h} b^{-1}
\underset{\Gamma_{\n}}{\neq} 1
$$
by Britton's Lemma (again because $|\n| \geqslant 2$).

\textbf{Case 3:} $\n \equiv k \ (\text{mod } \m_1^h d)$. \\
Let us write $\n - k = \ell \m_1^h d$. Then $w \underset{\Gamma_{\n}}{=} b^{\ell \n_1^h d} b b^{-\ell \n_1^h d}
b^{-1} \underset{\Gamma_{\n}}{=} 1$. \qed

\smallskip

\textbf{Proof of Theorem \ref{ConvSub}.} The "if" is a particular case of Theorem \ref{ConvSubOnlyIf}
below.

We prove now the "only if" part. Let $\Gamma_n = BS(\m,\k_n)$. We assume the sequence
$(\Gamma_n)_n$ to converge and condition (a) not to hold. We have to show that condition (b) holds.

Fix $h \geqslant 1$. For $k \in \Z$ we set
$$
w_k = a^{h+1} b^\m a^{-1} b^{-k} a^{-h} b a^{h+1} b^\m a^{-1} b^k a^{-h} b^{-1} \ .
$$
As $(\Gamma_n)_n$ converges, we have (for each $k \in \Z$) either $w_k \underset{\Gamma_n}{=} 1$ for $n$ large enough,
or $w_k \underset{\Gamma_n}{\neq} 1$ for $n$ large enough. As $\k_n$ is relatively prime to $\m$ for all $n$, Lemma
\ref{Congruence} ensures that (for each $k \in \Z$) either $\k_n \equiv k \ (\text{mod } \m^h)$ for $n$ large enough,
or $\k_n \not\equiv k \ (\text{mod } \m^h)$ for $n$ large enough. This implies that $\k_n$ is eventually constant
modulo $\m^h$ ($\forall h \geqslant 1$).

It remains to show that $|\k_n| \to \infty$. Assume by contradiction that there exists some $\ell \in \Z$ such
that $\k_n = \ell$ for infinitely many $n$. As (a) does not hold (i.e $(\k_n)_n$ is not eventually constant),
it is sufficient to treat the two following cases:

\textbf{Case 1:} $\exists \ell' \neq \ell$ such that $\k_n = \ell'$ for infinitely many $n$. \\
Take $h$ large enough so that $\m^h > |\ell - \ell'|$. The sequence $\k_n$ cannot be eventually constant
modulo $\m^h$, in contradiction with the first part of the proof.

\textbf{Case 2:} $\exists$ a subsequence $(\k_{n_j})_j$ of $(\k_n)_n$ such that $|\k_{n_j}| \to \infty$. \\
We set $w = ab^\m a^{-1}b^{-\ell}$. For infinitely many $n$ (those values for which $\k_n = \ell$), we have
$w \underset{\Gamma_n}{=} 1$. On the other hand $|\k_{n_j}| > \ell$ for $j$ large enough. For these values
of $j$, we have
$$
    w \underset{\Gamma_{n_j}}{=} b^{\k_{n_j} - \ell} \underset{\Gamma_{n_j}}{\neq} 1 \ .
$$
This contradicts the assumption on the sequence $(\Gamma_n)_n$ to converge. \qed

\smallskip

What remains now to do is to prove the following Theorem, which is a little bit more general than the
"if" part of Theorem \ref{ConvSub}. The proof will need some preliminary lemmas.

\begin{Thm}\label{ConvSubOnlyIf}
    Let $\m \in \Z^*$ and let $(\k_n)_n$ be a sequence in $\Z^*$. If $|\k_n| \to \infty$ and if $\forall h
    \geqslant 1$ the sequence $(\k_n)_n$ is eventually constant modulo $\m^h$, then the sequence $(BS(\m,\k_n))_n$ is
    convergent in ${\cal G}_2$.
\end{Thm}

\begin{Lem}\label{Prelim1}
    Let $\m,\n,\n' \in \Z^*$ and $h \geqslant 1$. If $\n \equiv \n'$ $(\text{mod } \m^h)$, there exists
    $s_0, \ldots, s_h$; $s'_0, \ldots, s'_h$; $r_1, \ldots, r_h$, which are unique, such that:
    \begin{enumerate}
        \item[(i)] $0 \leqslant r_i < \m \ \forall i$; $s_0 = 1 = s'_0$;
        \item[(ii)] $s_{i-1} \n = s_i \m + r_i $ and $s'_{i-1} \n' = s'_i \m + r_i$ $\forall \, 1\leqslant i \leqslant h$;
        \item[(iii)] $s_i \equiv s'_i \ (\text{mod } \m^{h-i}) \ \forall \, 0 \leqslant i \leqslant h$.
    \end{enumerate}

\end{Lem}

\textbf{Proof.} Given the congruence $\n \equiv \n'$ $(\text{mod } \m^h)$, we obtain (by Euclidean division)
$s_{0} \n = \n = s_1 \m + r_1 $ and $s'_{0} \n' = \n' = s'_1 \m + r_1$ with $0 \leqslant r_1 \leqslant \m$
and $s_1 \equiv s'_1$ $(\text{mod } \m^{h-1})$. Hence we have $s_1 \n \equiv s'_1 \n'$ $(\text{mod } \m^{h-1})$.
(Let us emphasize that we do not necessary have $s_1 \n \equiv s'_1 \n'$ $(\text{mod } \m^h)$.)

Now, it just remains to iterate the above and uniqueness follows from construction. \qed

\smallskip

Given a word $w$ in $\F_2$, we may use Britton's Lemma to reduce it in $BS(\m,\n)$ or $BS(\m,\n')$. But $w$
could be reducible in one of these groups and not in the other one. Even if it is reducible in both groups
the result is not the same word in general. The purpose of next statement is, under some assumptions, to
control the parallel process of reduction in both groups. This will be useful to ensure that $w$ is a
relation in $BS(\m,\n)$ if and only if it is one in $BS(\m,\n')$ (under some assumptions).

\begin{Lem}\label{Prelim2}
    Let $\m,\n,\n' \in \Z^*$ and $h \geqslant m \geqslant 1$. Assume that $\n \equiv \n'$ $(\text{mod } \m^h)$
    and let
    \begin{eqnarray*}
        \alpha & = & k_0 + k_1 \n + k_2 s_1 \n + \ldots + k_m s_{m-1} \n  \\
        \alpha' & = & k_0 + k_1 \n' + k_2 s'_1 \n' + \ldots + k_m s'_{m-1} \n'
    \end{eqnarray*}
    where $s_0, \ldots, s_h$; $s'_0, \ldots, s'_h$; $r_1, \ldots, r_h$ are given by Lemma \ref{Prelim1}. We
    assume moreover that we have $|k_0| < \min(|\n|,|\n'|)$.
    \begin{enumerate}
        \item[(i)] We have $\alpha \equiv 0$ $(\text{mod } \m)$ if and only if $\alpha' \equiv 0$ $(\text{mod } \m)$.
        If this happens we get $ab^\alpha a^{-1} \underset{BS(\m,\n)}{=} b^\beta$ and $ab^{\alpha'} a^{-1}
        \underset{BS(\m,\n')}{=} b^{\beta'}$ with
        \begin{eqnarray*}
            \beta & = & \ell_1 \n + \ell_2 s_1 \n + \ldots + \ell_{m+1} s_{m} \n  \\
            \beta' & = & \ell_1 \n' + \ell_2 s'_1 \n' + \ldots + \ell_{m+1} s'_{m} \n' \ .
        \end{eqnarray*}
        \item[(ii)] We have $\alpha \equiv 0$ $(\text{mod } \n)$ if and only if $\alpha' \equiv 0$ $(\text{mod } \n')$.
        If this happens we get $a^{-1}b^\alpha a \underset{BS(\m,\n)}{=} b^\beta$ and $a^{-1} b^{\alpha'} a
        \underset{BS(\m,\n')}{=} b^{\beta'}$ with
        \begin{eqnarray*}
            \beta & = & \ell_0 + \ell_1 \n + \ell_2 s_1 \n + \ldots + \ell_{m-1} s_{m-2} \n  \\
            \beta' & = & \ell_0 + \ell_1 \n' + \ell_2 s'_1 \n' + \ldots + \ell_{m-1} s'_{m-2} \n' \ .
        \end{eqnarray*}
    \end{enumerate}
\end{Lem}

\textbf{Proof.} (i) We have $\alpha \equiv \alpha' \ (\text{mod } \m)$ by construction. Assume now that
$\alpha \equiv 0$ and $ \alpha' \equiv 0 \ (\text{mod } \m)$. We have
$$
\alpha = k_0 + k_1 r_1 + \ldots + k_m r_m + k_1 s_1 \m + \ldots + k_m s_m \m \ .
$$
As $\alpha \equiv 0 \ (\text{mod } \m)$, we obtain $ab^\alpha a^{-1} \underset{BS(\m,\n)}{=} b^\beta$ with
$$
\beta = \frac{\n}{\m}(k_0 + k_1 r_1 + \ldots + k_m r_m) + k_1 s_1 \n + \ldots + k_m s_m \n
$$
Thus we set $\ell_1 = \frac{1}{\m}(k_0 + k_1 r_1 + \ldots + k_m r_m)$ and $\ell_i = k_{i-1}$ for $2 \leqslant
i \leqslant m+1$, and doing the same calculation with $\alpha'$ in $BS(\m,\n')$, we obtain also
$$
\beta' = \ell_1 \n' + \ell_2 s'_1 \n' + \ldots + \ell_{m+1} s'_{m} \n' \ .
$$

(ii) As $|\n|> |k_0|$ and $|\n'| > |k_0|$, we have $\alpha \equiv 0 \ (\text{mod } \n)$ if and only if $k_0 = 0$
if and only if $\alpha' \equiv 0 \ (\text{mod } \n')$. Suppose now that it is the case. We have $ab^\alpha
a^{-1} = b^\beta$ in $BS(\m,\n)$ with
\begin{eqnarray*}
\beta & = & k_1 \m + k_2 s_1 \m + \ldots + k_m s_{m-1} \m  \\
      & = & k_1 \m - k_2 r_1 - \ldots - k_m r_{m-1} + k_2 \n + k_3 s_1 \n + \ldots + k_m s_{m-2} \n \ .
\end{eqnarray*}
Hence we set $\ell_0 = k_1 \m - k_2 r_1 - \ldots - k_m r_{m-1}$ and $\ell_i = k_{i+1}$ for $1 \leqslant
i\leqslant m-1$. Again, doing the same calculation with $\alpha'$ in $BS(\m,\n')$, we obtain also
$$
\beta' = \ell_0 + \ell_1 \n' + \ell_2 s'_1 \n' + \ldots + \ell_{m-1} s'_{m-2} \n' \ .
$$
This completes the proof. \qed

\begin{Lem}\label{Prelim3}
    Let $\m \in \Z^*$, let $(\k_n)_n$ be a sequence in $\Z^*$ such that $|\k_n| \to \infty$ and $\forall h
    \geqslant 1$ $(\k_n)_n$ is eventually constant modulo $\m^h$ and let $w \in \F_2$. Then, the following alternative
    holds:
    \begin{enumerate}
        \item[(a)] either $w = b^{\lambda_n}$ in $BS(\m,\k_n)$ for $n$ large enough;
        \item[(b)] or $w$ is in $BS(\m,\k_n) \setminus \langle b \rangle$ for $n$ large enough.
    \end{enumerate}
\end{Lem}

\textbf{Proof.} We define $\Gamma_n = BS(\m,\k_n)$. Let us write $w = b^{\alpha_0} a^{\varepsilon_1} b^{\alpha_1}
\cdots a^{\varepsilon_m} b^{\alpha_m}$ with $\varepsilon_i = \pm 1$ and $\alpha_i \in \Z$, reduced in
the sense that $\alpha_i = 0$ implies $\varepsilon_{i+1} = \varepsilon_{i}$ for all $i \in \{1, \ldots, m-1 \}$.
We assume (b) not to hold, i.e. $w= b^{\lambda_n}$ in $\Gamma_n$ for infinitely many $n$. Then the sum
$\varepsilon_1 + \ldots + \varepsilon_m$ has clearly to be zero (in particular $m$ is even). We have to show
that $w= b^{\lambda_n}$ for $n$ large enough.

For $n$ large enough, we may assume that, $|\k_n| > |\alpha_j|$  for all $1 \leqslant j \leqslant m$
and the $\k_n$'s are all congruent modulo $\m^m$. We take a value of $n$ such that moreover $w= b^{\lambda_n}$ in
$\Gamma_n$ (there are infinitely many ones) and apply Britton's Lemma. This ensures the existence of
an index $j$ such that $\varepsilon_j = 1 = -\varepsilon_{j+1}$ and $\alpha_j \equiv 0 \ (\text{mod } \m)$
(since $|\k_n| > |\alpha_j|$ for all $j$). By Lemma \ref{Prelim2}, for all $n$ large enough
$$
w \underset{\Gamma_n}{=} b^{\alpha_0} \cdots a^{\varepsilon_{j-1}} b^{\alpha_{j-1} + \beta_j + \alpha_{j+1}}
a^{\varepsilon_{j+2}} \cdots b^{\alpha_m}
$$
with $\beta_j = \ell_1 \k_n$ (depending on $n$). Hence we are allowed to write
$$
w \underset{\Gamma_n}{=} b^{\alpha'_{0,n}} a^{\varepsilon'_1} b^{\alpha'_{1,n}}
\cdots a^{\varepsilon'_{m-2,n}} b^{\alpha'_{m-2,n}}
$$
for $n$ large enough, with $\varepsilon'_i = \pm 1$ and $\alpha'_{j,n} = k'_{0,j} + k'_{1,j} \k_n$, where the
$k'_{i,j}$'s do not depend on $n$.

Now, for $n$ large enough, we may assume that, $|\k_n| > |k'_{0,j}|$ for all $1 \leqslant j \leqslant m-1$
(and the $\k_n$'s are all congruent modulo $\m^m$). Again we take a value of $n$ such that moreover $w=b^{\lambda_n}$ in
$\Gamma_n$ and apply Britton's Lemma. This ensures the existence of an index $j$ such that either
$\varepsilon'_j = 1 = -\varepsilon'_{j+1}$ and $\alpha'_{j,n} \equiv 0 \ (\text{mod } \m)$, or
$\varepsilon'_j = -1 = -\varepsilon'_{j+1}$ and $\alpha'_{j,n} \equiv 0 \ (\text{mod } \k_n)$. In both cases,
while applying Lemma \ref{Prelim2}, we obtain
$$
w \underset{\Gamma_n}{=} b^{\alpha''_{0,n}} a^{\varepsilon''_1} b^{\alpha''_{1,n}}
\cdots a^{\varepsilon''_{m-4,n}} b^{\alpha''_{m-4,n}}
$$
for $n$ large enough, with $\varepsilon''_i = \pm 1$ and $\alpha''_{j,n} = k''_{0,j} + k''_{1,j} \k_n + k''_{2,j} s_{1,n} \k_n$,
where the $k''_{i,j}$'s do not depend on $n$.

And so on, and so forth, setting $m' = \frac{m}{2}$, we get finally $w = b^{\alpha^{(m')}_{0,n}}$ in $\Gamma_n$
for $n$ large enough, with
$$
\alpha^{(m')}_{0,n} = k^{(m')}_{0,0} + k^{(m')}_{1,0} \k_n + k^{(m')}_{2,0} s_{1,n} \k_n +
\ldots + k^{(m')}_{m',0} s_{m'-1,n} \k_n
$$
where the $k^{m'}_{i,0}$'s do not depend on $n$. It only remains to set $\lambda_n = \alpha^{(m')}_{0,n}$.
\qed

\smallskip

Let us now introduce the homomorphisms $\psi_\n : BS(\m,\n) \to \text{Aff}(\R)$ (for $\n \in \Z^*$) given by
$\psi_\n (a)(x) = \frac{\n}{\m} x$ and $\psi_\n (b)(x) = x+1$.

\begin{Lem}\label{Prelim4}
    Let $w \in \F_2$. We have either $\psi_\n (w) = 1$ for $|\n|$ large enough or $\psi_\n (w) \neq 1$
    for $|\n|$ large enough.
\end{Lem}

\textbf{Proof.} Let us write $w = b^{\alpha_0} a^{\varepsilon_1} b^{\alpha_1} \ldots a^{\varepsilon_k}
b^{\alpha_k}$ with $\varepsilon_i = \pm 1$ and $\alpha_i \in \Z$. Set next $\sigma_0 = 0$, $\sigma_i =
\varepsilon_1 + \ldots + \varepsilon_i$ for $1 \leqslant i \leqslant k$ and $m = \min_{0 \leqslant i
\leqslant k} \sigma_i$. We get by calculation that
$$
\psi_\n (w)(x) = \left( \frac{\n}{\m} \right)^{\sigma_k} x + \left( \frac{\n}{\m} \right)^m P_w \left(\frac{\n}{\m} \right)
$$
where $P_w$ is the polynomial defined by $P_w (y) = \sum_{i=0}^k \alpha_i y^{\sigma_i - m}$. Let us assume the
second term of alternative not to hold, i.e. $\psi_\n(w) = 1$ for infinitely many values of $\n$. Hence we have
$\sigma_k = 0$ and for all those values of $\n$, $P_w(\frac{\n}{\m}) = 0$. As $P_w$ is a polynomial with infinitely many
roots, it is the zero polynomial. This shows that $\psi_\n(w) = 1$ for all $\n$. \qed

\smallskip

\textbf{Proof of Theorem \ref{ConvSubOnlyIf}.} It is easy to show that a word $w$ is equal to $1$ in
$BS(\m,\n)$ if and only if it is in the subgroup generated by $b$ and  $\psi_\n(w) = 1$. It is also
a consequence of (the proof of) Theorem 1 in \cite{GalJan}.
Let $w\in \F_2$. Lemmas \ref{Prelim3} and \ref{Prelim4} immediately imply that either $w=1$ in $BS(\m,\k_n)$
for $n$ large enough or $w \neq 1$ in $BS(\m,\k_n)$ for $n$ large enough. \qed

Theorem \ref{ConvSub} is now completely established.

\section{Baumslag-Solitar groups and $\m$-adic integers}\label{Topol}

This last section is entirely devoted to the proof of Theorem \ref{UnifCont}. We show first that the map
$\Psi$ is uniformly continuous. In view of the distances we put on ${\cal G}_2$ and $\Z_\m^\times$, it is equivalent
to show that for any $h \geqslant 1$ there exists $r \geqslant 1$ such that we have $\n \equiv \n' \ (\text{mod } \m^h)$
whenever $BS(\m,\n)$ and $BS(\m,\n')$ (taken in $X_\m$) have the same relations up to length $r$.

Fix $h\geqslant 1$. Our candidate is $r = 2\m^h + 4h + 2\m + 4$. Assume that $BS(\m,\n)$ and $BS(\m,\n')$ have the same
relations up to length $r$. For $0\leqslant k \leqslant \m^h-1$ let
$$
w_k = a^{h+1} b^\m a^{-1} b^{-k} a^{-h} b a^{h+1} b^\m a^{-1} b^k a^{-h} b^{-1} \ .
$$
(Remark that these words are exactly those which appear in the proof of the "only if" part of theorem \ref{ConvSub}.
We are in fact improving this proof in order to get the uniform continuity.) We have $\ell(w_k) \leqslant r$
for all $k$. Having by assumption $w = 1$ in $BS(\m,\n)$ if and only if $w = 1$ in $BS(\m,\n')$, Lemma
\ref{Congruence} implies $\n \equiv \n' \ (\text{mod } \m^h)$.

The space $\Z_\m^\times$ being complete and the uniform continuity of $\Psi$ being now proved, the existence
of a (unique) uniformly continuous extension $\overline{\Psi}$ on $\overline{X_\m}$ is a standard fact
(see \cite{Dug}, Chapter XIV, Theorem 5.2. for instance).

Let us now show that $\overline{\Psi}$ is surjective. The space $\overline{X_\m}$ being compact,
$\text{im}(\overline{\Psi})$ is closed in $\Z_\m^\times$. Moreover it is dense since it contains the set of
$\n$'s relatively prime to $\m$ and such that $|\n| \geqslant 2$.

Finally, we consider the sequence $(BS(\m,1+\m+\m^n))_n$, which is convergent by theorem \ref{ConvSub} and we
call the limit $\Gamma$. We have $\overline{\Psi}(\Gamma) = 1+\m = \overline{\Psi}(BS(\m,1+\m))$. On the other
hand, we have $\Gamma \neq BS(1+\m)$, since $ab^\m a^{-1}b^{-(\m+1)} = 1$ in $BS(\m,1+\m)$ while
$ab^\m a^{-1}b^{-(\m+1)} \neq 1$ in $BS(\m,1+\m+\m^n)$ for all $n$. This is the non-injectivity of
$\overline{\Psi}$ and completes the proof. \qed


\bibliographystyle{alpha}

\medskip

Author address:

\vspace{2mm}

Institut de Math\'ematiques - Universit\'e de Neuch\^atel

Rue Emile Argand 11

CH-2007 Neuch\^atel - SWITZERLAND

\vspace{2mm}

yves.stalder@unine.ch

Phone: +41 32 718 28 17; Fax: +41 32 718 28 01

\end{document}